\documentclass{article}
\usepackage{amssymb,amsmath,theorem,euscript}

\usepackage{upgreek}
\usepackage{calrsfs}

\input{epsf}

\newcounter{sec}

\def\sm{\smallskip}

\newcounter{punct}[sec]

\def\punct{\refstepcounter{punct}{\arabic{sec}.\arabic{punct}.  }}

\def\COUNTERS{\addtocounter{sec}{1}
              \setcounter{punct}{0}
          \setcounter{equation}{0}
          \setcounter{theorem}{0}
                  }
\newtheorem{theorem}{Theorem}[sec]

\newtheorem{lemma}[theorem]{Lemma}

\begin{document}

 \def\ov{\overline}
\def\wt{\widetilde}
\def\wh{\widehat}
 \newcommand{\sgn}{\mathop {\mathrm {sign}}\nolimits}
 \newcommand{\rk}{\mathop {\mathrm {rk}}\nolimits}
\newcommand{\Aut}{\mathop {\mathrm {Aut}}\nolimits}
\newcommand{\Out}{\mathop {\mathrm {Out}}\nolimits}
\newcommand{\Abs}{\mathop {\mathrm {Abs}}\nolimits}
\renewcommand{\Re}{\mathop {\mathrm {Re}}\nolimits}
\renewcommand{\Im}{\mathop {\mathrm {Im}}\nolimits}
 \newcommand{\tr}{\mathop {\mathrm {tr}}\nolimits}
  \newcommand{\Hom}{\mathop {\mathrm {Hom}}\nolimits}
   \newcommand{\diag}{\mathop {\mathrm {diag}}\nolimits}
   \newcommand{\supp}{\mathop {\mathrm {supp}}\nolimits}
 \newcommand{\im}{\mathop {\mathrm {im}}\nolimits}
 \newcommand{\grad}{\mathop {\mathrm {grad}}\nolimits}
  \newcommand{\sgrad}{\mathop {\mathrm {sgrad}}\nolimits}
 \newcommand{\rot}{\mathop {\mathrm {rot}}\nolimits}
  \renewcommand{\div}{\mathop {\mathrm {div}}\nolimits}

\def\Br{\mathrm {Br}}
\def\Vir{\mathrm {Vir}}

 \def\Ham{\mathrm {Ham}}
\def\SL{\mathrm {SL}}
\def\Pol{\mathrm {Pol}}
\def\SU{\mathrm {SU}}
\def\GL{\mathrm {GL}}
\def\U{\mathrm U}
\def\OO{\mathrm O}
 \def\Sp{\mathrm {Sp}}
  \def\Ad{\mathrm {Ad}}
 \def\SO{\mathrm {SO}}
\def\SOS{\mathrm {SO}^*}
 \def\Diff{\mathrm{Diff}}
 \def\Vect{\mathfrak{Vect}}
\def\PGL{\mathrm {PGL}}
\def\PU{\mathrm {PU}}
\def\PSL{\mathrm {PSL}}
\def\Symp{\mathrm{Symp}}
\def\Cont{\mathrm{Cont}}
\def\End{\mathrm{End}}
\def\Mor{\mathrm{Mor}}
\def\Aut{\mathrm{Aut}}
 \def\PB{\mathrm{PB}}
\def\Fl{\mathrm {Fl}}
\def\Symm{\mathrm {Symm}} 
 \def\Herm{\mathrm {Herm}} 
  \def\SDiff{\mathrm {SDiff}} 
 
 \def\cA{\mathcal A}
\def\cB{\mathcal B}
\def\cC{\mathcal C}
\def\cD{\mathcal D}
\def\cE{\mathcal E}
\def\cF{\mathcal F}
\def\cG{\mathcal G}
\def\cH{\mathcal H}
\def\cJ{\mathcal J}
\def\cI{\mathcal I}
\def\cK{\mathcal K}
 \def\cL{\mathcal L}
\def\cM{\mathcal M}
\def\cN{\mathcal N}
 \def\cO{\mathcal O}
\def\cP{\mathcal P}
\def\cQ{\mathcal Q}
\def\cR{\mathcal R}
\def\cS{\mathcal S}
\def\cT{\mathcal T}
\def\cU{\mathcal U}
\def\cV{\mathcal V}
 \def\cW{\mathcal W}
\def\cX{\mathcal X}
 \def\cY{\mathcal Y}
 \def\cZ{\mathcal Z}
\def\0{{\ov 0}}
 \def\1{{\ov 1}}
 \def\frA{\mathfrak A}
 \def\frB{\mathfrak B}
\def\frC{\mathfrak C}
\def\frD{\mathfrak D}
\def\frE{\mathfrak E}
\def\frF{\mathfrak F}
\def\frG{\mathfrak G}
\def\frH{\mathfrak H}
\def\frI{\mathfrak I}
 \def\frJ{\mathfrak J}
 \def\frK{\mathfrak K}
 \def\frL{\mathfrak L}
\def\frM{\mathfrak M}
 \def\frN{\mathfrak N} \def\frO{\mathfrak O} \def\frP{\mathfrak P} \def\frQ{\mathfrak Q} \def\frR{\mathfrak R}
 \def\frS{\mathfrak S} \def\frT{\mathfrak T} \def\frU{\mathfrak U} \def\frV{\mathfrak V} \def\frW{\mathfrak W}
 \def\frX{\mathfrak X} \def\frY{\mathfrak Y} \def\frZ{\mathfrak Z} \def\fra{\mathfrak a} \def\frb{\mathfrak b}
 \def\frc{\mathfrak c} \def\frd{\mathfrak d} \def\fre{\mathfrak e} \def\frf{\mathfrak f} \def\frg{\mathfrak g}
 \def\frh{\mathfrak h} \def\fri{\mathfrak i} \def\frj{\mathfrak j} \def\frk{\mathfrak k} \def\frl{\mathfrak l}
 \def\frm{\mathfrak m} \def\frn{\mathfrak n} \def\fro{\mathfrak o} \def\frp{\mathfrak p} \def\frq{\mathfrak q}
 \def\frr{\mathfrak r} \def\frs{\mathfrak s} \def\frt{\mathfrak t} \def\fru{\mathfrak u} \def\frv{\mathfrak v}
 \def\frw{\mathfrak w} \def\frx{\mathfrak x} \def\fry{\mathfrak y} \def\frz{\mathfrak z} \def\frsp{\mathfrak{sp}}
 \def\bfa{\mathbf a} \def\bfb{\mathbf b} \def\bfc{\mathbf c} \def\bfd{\mathbf d} \def\bfe{\mathbf e} \def\bff{\mathbf f}
 \def\bfg{\mathbf g} \def\bfh{\mathbf h} \def\bfi{\mathbf i} \def\bfj{\mathbf j} \def\bfk{\mathbf k} \def\bfl{\mathbf l}
 \def\bfm{\mathbf m} \def\bfn{\mathbf n} \def\bfo{\mathbf o} \def\bfp{\mathbf p} \def\bfq{\mathbf q} \def\bfr{\mathbf r}
 \def\bfs{\mathbf s} \def\bft{\mathbf t} \def\bfu{\mathbf u} \def\bfv{\mathbf v} \def\bfw{\mathbf w} \def\bfx{\mathbf x}
 \def\bfy{\mathbf y} \def\bfz{\mathbf z} \def\bfA{\mathbf A} \def\bfB{\mathbf B} \def\bfC{\mathbf C} \def\bfD{\mathbf D}
 \def\bfE{\mathbf E} \def\bfF{\mathbf F} \def\bfG{\mathbf G} \def\bfH{\mathbf H} \def\bfI{\mathbf I} \def\bfJ{\mathbf J}
 \def\bfK{\mathbf K} \def\bfL{\mathbf L} \def\bfM{\mathbf M} \def\bfN{\mathbf N} \def\bfO{\mathbf O} \def\bfP{\mathbf P}
 \def\bfQ{\mathbf Q} \def\bfR{\mathbf R} \def\bfS{\mathbf S} \def\bfT{\mathbf T} \def\bfU{\mathbf U} \def\bfV{\mathbf V}
 \def\bfW{\mathbf W} \def\bfX{\mathbf X} \def\bfY{\mathbf Y} \def\bfZ{\mathbf Z} \def\bfw{\mathbf w}
 \def\R {{\mathbb R }} \def\C {{\mathbb C }} \def\Z{{\mathbb Z}} \def\H{{\mathbb H}} \def\K{{\mathbb K}}
 \def\N{{\mathbb N}} \def\Q{{\mathbb Q}} \def\A{{\mathbb A}} \def\T{\mathbb T} \def\P{\mathbb P} \def\G{\mathbb G}
 \def\bbA{\mathbb A} \def\bbB{\mathbb B} \def\bbD{\mathbb D} \def\bbE{\mathbb E} \def\bbF{\mathbb F} \def\bbG{\mathbb G}
 \def\bbI{\mathbb I} \def\bbJ{\mathbb J} \def\bbL{\mathbb L} \def\bbM{\mathbb M} \def\bbN{\mathbb N} \def\bbO{\mathbb O}
 \def\bbP{\mathbb P} \def\bbQ{\mathbb Q} \def\bbS{\mathbb S} \def\bbT{\mathbb T} \def\bbU{\mathbb U} \def\bbV{\mathbb V}
 \def\bbW{\mathbb W} \def\bbX{\mathbb X} \def\bbY{\mathbb Y} \def\kappa{\varkappa} \def\epsilon{\varepsilon}
 \def\phi{\varphi} \def\le{\leqslant} \def\ge{\geqslant}

\def\UU{\bbU}
\def\Mat{\mathrm{Mat}}
\def\tto{\rightrightarrows}

\def\F{\mathbf{F}}

\def\Gms{\mathrm {Gms}}
\def\Ams{\mathrm {Ams}}
\def\Isom{\mathrm {Isom}}

\def\Gr{\mathrm{Gr}}

\def\graph{\mathrm{graph}}

\def\la{\langle}
\def\ra{\rangle}


 \def\ov{\overline}
\def\wt{\widetilde}

\renewcommand{\Re}{\mathop {\mathrm {Re}}\nolimits}
\def\Br{\mathrm {Br}}

  \def\Match{\mathrm {Match}}
 \def\Isom{\mathrm {Isom}}
 \def\Hier{\mathrm {Hier}}
\def\SL{\mathrm {SL}}
\def\SU{\mathrm {SU}}
\def\GL{\mathrm {GL}}
\def\U{\mathrm U}
\def\OO{\mathrm O}
 \def\Sp{\mathrm {Sp}}
  \def\GLO{\mathrm {GLO}}
 \def\SO{\mathrm {SO}}
\def\SOS{\mathrm {SO}^*}
 \def\Diff{\mathrm{Diff}}
 \def\Vect{\mathfrak{Vect}}
\def\PGL{\mathrm {PGL}}
\def\PU{\mathrm {PU}}
\def\PSL{\mathrm {PSL}}
\def\Symp{\mathrm{Symp}}
\def\ASymm{\mathrm{Asymm}}
\def\Asymm{\mathrm{Asymm}}
\def\Gal{\mathrm{Gal}}
\def\End{\mathrm{End}}
\def\Mor{\mathrm{Mor}}
\def\Aut{\mathrm{Aut}}
 \def\PB{\mathrm{PB}}
 \def\cA{\mathcal A}
\def\cB{\mathcal B}
\def\cC{\mathcal C}
\def\cD{\mathcal D}
\def\cE{\mathcal E}
\def\cF{\mathcal F}
\def\cG{\mathcal G}
\def\cH{\mathcal H}
\def\cJ{\mathcal J}
\def\cI{\mathcal I}
\def\cK{\mathcal K}
 \def\cL{\mathcal L}
\def\cM{\mathcal M}
\def\cN{\mathcal N}
 \def\cO{\mathcal O}
\def\cP{\mathcal P}
\def\cQ{\mathcal Q}
\def\cR{\mathcal R}
\def\cS{\mathcal S}
\def\cT{\mathcal T}
\def\cU{\mathcal U}
\def\cV{\mathcal V}
 \def\cW{\mathcal W}
\def\cX{\mathcal X}
 \def\cY{\mathcal Y}
 \def\cZ{\mathcal Z}
\def\0{{\ov 0}}
 \def\1{{\ov 1}}
 
 \def\frA{\mathfrak A}
 \def\frB{\mathfrak B}
\def\frC{\mathfrak C}
\def\frD{\mathfrak D}
\def\frE{\mathfrak E}
\def\frF{\mathfrak F}
\def\frG{\mathfrak G}
\def\frH{\mathfrak H}
\def\frI{\mathfrak I}
 \def\frJ{\mathfrak J}
 \def\frK{\mathfrak K}
 \def\frL{\mathfrak L}
\def\frM{\mathfrak M}
 \def\frN{\mathfrak N} \def\frO{\mathfrak O} \def\frP{\mathfrak P} \def\frQ{\mathfrak Q} \def\frR{\mathfrak R}
 \def\frS{\mathfrak S} \def\frT{\mathfrak T} \def\frU{\mathfrak U} \def\frV{\mathfrak V} \def\frW{\mathfrak W}
 \def\frX{\mathfrak X} \def\frY{\mathfrak Y} \def\frZ{\mathfrak Z} \def\fra{\mathfrak a} \def\frb{\mathfrak b}
 \def\frc{\mathfrak c} \def\frd{\mathfrak d} \def\fre{\mathfrak e} \def\frf{\mathfrak f} \def\frg{\mathfrak g}
 \def\frh{\mathfrak h} \def\fri{\mathfrak i} \def\frj{\mathfrak j} \def\frk{\mathfrak k} \def\frl{\mathfrak l}
 \def\frm{\mathfrak m} \def\frn{\mathfrak n} \def\fro{\mathfrak o} \def\frp{\mathfrak p} \def\frq{\mathfrak q}
 \def\frr{\mathfrak r} \def\frs{\mathfrak s} \def\frt{\mathfrak t} \def\fru{\mathfrak u} \def\frv{\mathfrak v}
 \def\frw{\mathfrak w} \def\frx{\mathfrak x} \def\fry{\mathfrak y} \def\frz{\mathfrak z} \def\frsp{\mathfrak{sp}}
 \def\bfa{\mathbf a} \def\bfb{\mathbf b} \def\bfc{\mathbf c} \def\bfd{\mathbf d} \def\bfe{\mathbf e} \def\bff{\mathbf f}
 \def\bfg{\mathbf g} \def\bfh{\mathbf h} \def\bfi{\mathbf i} \def\bfj{\mathbf j} \def\bfk{\mathbf k} \def\bfl{\mathbf l}
 \def\bfm{\mathbf m} \def\bfn{\mathbf n} \def\bfo{\mathbf o} \def\bfp{\mathbf p} \def\bfq{\mathbf q} \def\bfr{\mathbf r}
 \def\bfs{\mathbf s} \def\bft{\mathbf t} \def\bfu{\mathbf u} \def\bfv{\mathbf v} \def\bfw{\mathbf w} \def\bfx{\mathbf x}
 \def\bfy{\mathbf y} \def\bfz{\mathbf z} \def\bfA{\mathbf A} \def\bfB{\mathbf B} \def\bfC{\mathbf C} \def\bfD{\mathbf D}
 \def\bfE{\mathbf E} \def\bfF{\mathbf F} \def\bfG{\mathbf G} \def\bfH{\mathbf H} \def\bfI{\mathbf I} \def\bfJ{\mathbf J}
 \def\bfK{\mathbf K} \def\bfL{\mathbf L} \def\bfM{\mathbf M} \def\bfN{\mathbf N} \def\bfO{\mathbf O} \def\bfP{\mathbf P}
 \def\bfQ{\mathbf Q} \def\bfR{\mathbf R} \def\bfS{\mathbf S} \def\bfT{\mathbf T} \def\bfU{\mathbf U} \def\bfV{\mathbf V}
 \def\bfW{\mathbf W} \def\bfX{\mathbf X} \def\bfY{\mathbf Y} \def\bfZ{\mathbf Z} \def\bfw{\mathbf w}

 \def\R {{\mathbb R }} \def\C {{\mathbb C }} \def\Z{{\mathbb Z}} \def\H{{\mathbb H}} \def\K{{\mathbb K}}
 \def\N{{\mathbb N}} \def\Q{{\mathbb Q}} \def\A{{\mathbb A}} \def\T{\mathbb T} \def\P{\mathbb P} \def\G{\mathbb G}
 \def\bbA{\mathbb A} \def\bbB{\mathbb B} \def\bbD{\mathbb D} \def\bbE{\mathbb E} \def\bbF{\mathbb F} \def\bbG{\mathbb G}
 \def\bbI{\mathbb I} \def\bbJ{\mathbb J} \def\bbL{\mathbb L} \def\bbM{\mathbb M} \def\bbN{\mathbb N} \def\bbO{\mathbb O}
 \def\bbP{\mathbb P} \def\bbQ{\mathbb Q} \def\bbS{\mathbb S} \def\bbT{\mathbb T} \def\bbU{\mathbb U} \def\bbV{\mathbb V}
 \def\bbW{\mathbb W} \def\bbX{\mathbb X} \def\bbY{\mathbb Y} \def\kappa{\varkappa} \def\epsilon{\varepsilon}
 \def\phi{\varphi} \def\le{\leqslant} \def\ge{\geqslant}

\def\UU{\bbU}
\def\Mat{\mathrm{Mat}}
\def\tto{\rightrightarrows}

\def\Gr{\mathrm{Gr}}

\def\ch{\cosh}
\def\sh{\sinh}

\def\B{\bfB} 

\def\graph{\mathrm{graph}}

\def\gl{\mathfrak{gl}}

\def\la{\langle}
\def\ra{\rangle}

\def\V{\Updelta}
\def\ctg{\cot}

\def\e{\mathcal{E}}

\begin{center}
	\bf\Large
 Projectors separating  spectra
 \\
  for $L^2$ on symmetric spaces $\GL(n,\C)/\GL(n,\R)$
  
  \medspace
  
  \sc Yury A. Neretin%
  \footnote{Supported by the grants FWF, P25142, P28421}
  
\end{center}

{\small The Plancherel decomposition of $L^2$ on a pseudo-Riemannian
	symmetric space
	\newline $\GL(n,\C)/\GL(n,\R)$ has spectrum of  $[n/2]$ types. We write explicitly orthogonal projectors
	separating spectrum into uniform pieces.}

\section{Formulas for the projectors}

\COUNTERS

{\bf\punct Problem.} It is well-known that in various problems of non-commutative harmonic analysis
spectra split  into pieces of different nature. 

In 1977 I.~M.~Gelfand and S.~G.~Gindikin \cite{GG}
formulated a problem
 about an explicit description of a decomposition of $L^2$ on a semisimple Lie group $G$
  (and, more generally, on a semi-simple pseudo-Riemannian
symmetric space) 
into a direct sum of 
representations having uniform spectra. They also gave an answer for $G=\SL(2,\R)$ 
  in the terms of boundary values 
of holomorphic functions. After this work the question became a topic of intensive 
discussions and series of works.

The picture is quite nice for $L^2\bigl(\SL(2,\R)\bigr)$ and $L^2\bigl(\SL(2,\R)/H\bigr)$,
where $H$ is the diagonal subgroup, see \cite{Gin2}, \cite{Gin3},  \cite{AU}, and \cite{Mol1}, \cite{Ner-rest}
respectively.
 However, up to now the situation for higher groups  is far to be well-understood.

\sm

{\bf\punct Known approaches.} a) G.~I.~Olshanski \cite{Ols1} proposed a way to split off
holomorphic discrete series using boundary values of holomorphic functions, 
this approach was used in several works,
see, e.g. \cite{Kof}, \cite{KM}.

S.~G.~Gindikin, B.~Kr\o tz  and G.~\'Olafsson \cite{GK} showed that in some cases
an integral of the most continuous series
can be splitted off in a similar way.

\sm

b) V.~F.Molchanov \cite{Mol2} and S.~G.~Gindikin \cite{Gin3}
in different way obtained the desired decomposition for $L^2$ on 
multi-dimensional hyperboloids $\OO(p,q)/\OO(p,q-1)$
(this includes the cases $L^2\bigl(\SL(2,\R)\bigr)$ and $L^2\bigl(\SL(2,\R)/H\bigr)$ mentioned above).

\sm

c) In \cite{Ner-compl} there was proposed a way for separation of summands of complementary series, see more in
\cite{Ner-separ}, \cite{NO}.

\sm

{\bf \punct Purposes of the present work.}
For a pseudo-Riemannian symmetric space $\GL(n,\C)/\GL(n,\R)$
consider a decomposition of $L^2$ into a direct sum of subspaces,
in which $\GL(n,\C)$ has uniform spectra%
\footnote{Here $[n/2]$ denotes the integer part of $n/2$.},
\begin{equation}
L^2\bigl(\GL(n,\C)/\GL(n,\R)\bigr)=L_0\oplus L_1\oplus\dots\oplus L_{[n/2]}.
\label{eq:LLL}
\end{equation}
We also have a corresponding decomposition of the identity operator:
\begin{equation}
E=\Pi_0\oplus \Pi_1\oplus \dots \oplus \Pi_{[n/2]},
\label{eq:PPP}
\end{equation}
where $\Pi_r$ is the orthogonal projector to a subspace $L_r$.
 We intend to obtain  explicit expressions
for the  projectors $\Pi_r$.

The present work is  based on both results and auxiliary calculations of
Shigeru Sano \cite{San}.
The formulas for the projectors are given by Theorem \ref{th:Lambda}.

\sm

{\bf \punct The spaces $\GL(n,\C)/\GL(n,\R)$.} Below 
$$
G:=\GL(n,\C),\qquad H:=\GL(n,\R).
$$
We realize the symmetric space $G/H$ as the space $M$ of
all matrices $x\in\GL(n,\C)$ satisfying the condition
$$
x\ov x=1.
$$
The group $G$ acts on $M$ by  transformations
$$
g:\, x\mapsto x\circ g:=g^{-1} x \ov g,
$$
the stabilizer of $x=e$ is $H$. Notice that $H$ acts on $M$ by conjugations.

Denote by $C^\infty_c(G/H)$ the space of $C^\infty$-smooth functions on $G/H$ with compact support.
Let $\chi$ be an $H$-invariant distribution on $G/H$. It determines a 
$G$-intertwining operator 
$$A[\chi]: C^\infty_c(G/H)\to C^\infty(G/H)
$$
by the pairing%
\footnote{By $\la f, \chi\ra_L$ we denote a pairing of a test function
$f$ and a distribution $\chi$  on a manifold $L$.}
$$
A[\chi] f(g)=\la f(x\circ g,\chi\ra_{G/H},
$$
a function on $G$ obtained in this way is $H$-invariant and therefore it is a function on $G/H$.

Our purpose is to write  $H$-invariant distributions on $G/H$
determining the projectors $\Pi_r$.

\sm

{\bf \punct Notation.}
By $dx$ we denote a $G$-invariant measure on $G/H$ (it is unique up to a constant factor).
Denote by $S_m$ the symmetric group of order $m$,
by $\T^p$ the torus $(\R/2\pi \Z)^p$.
Denote by $\V$ the {\it Vandermonde expression}:
\begin{equation}
\V(s)=\V(s_1,\dots,s_n):=\prod_{1\le p<q\le n} (s_p-s_q).
\label{eq:vandermonde}
\end{equation}

\sm

{\bf\punct  Cartan subspaces and the Weyl integration formula.}
Let $\phi$, $\theta\in \R/2\pi \Z$, $t\in\R$. Denote by $v(\theta,t)$ 
a $2\times 2$ matrix given by 
$$
v(\theta,t)=
 e^{i\theta}
\begin{pmatrix}
\ch t&i\sh t\\-i \sh t&\ch t
\end{pmatrix}=
\begin{pmatrix}
 e^{i\theta} \ch t&i  e^{i\theta}\sh t\\-i  e^{i\theta} \sh t&  e^{i\theta}\ch t
\end{pmatrix}
,$$
its eigenvalues are
$$
e^{z}:=e^{t+i\theta}, \qquad e^{-\ov z}:=e^{-t+i\theta}.
$$

We define a {\it Cartan subspace} $A_k$, where $k=1$, 2, \dots, $[n/2]$, as the set of
 matrices $a\in M$ having the following block-diagonal form
\begin{equation}
a^{(k)}:=\begin{pmatrix}
    e^{i\phi_1}\\
    &\ddots\\
    && e^{i\phi_{n-2k}}\\
    &&& v(\theta_k,t_k)\\
     &&&&\ddots \\
      &&&&& v(\theta_1,t_1)
   \end{pmatrix}
   \label{eq:A-k}
\end{equation}
We equip $A_k$ with the standard Lebesgue measure
$$
da^{(k)}:=\prod_{l=1}^{n-2k} d\phi_l \prod_{m=1}^k (d\theta_m\, dt_m).
$$

The relative {\it Weyl group} $W_k$  corresponding to the Cartan subspace $A_k$ is
$$
W_k\simeq S_{n-2k} \times (S_k\ltimes \Z_2^k).
$$
The factor $S_{n-2k}$ acts on $A_k$ by permutations of $e^{i\phi_2}$, \dots, $e^{i\phi_{n-2k}}$,
the group
$S_k$ acts by permutations 
of pairs $(\theta_1,t_1)$, \dots, $(\theta_k,t_k)$.
The group $\Z_2^k$
is generated by reflections
\begin{equation}
 R_m:(t_1,\dots,t_{m-1},t_m,t_{m+1},\dots, t_k)\mapsto (t_1,\dots,t_{m-1},-t_m,t_{m+1},\dots, t_k)
 \label{eq:Rm}
\end{equation}
(all the coordinates $\phi$ and $\theta$ remain to be fixed).
We denote
\begin{equation}
\gamma_k:=\frac1{\# W_k}=\frac 1{k! (n-2k)!\, 2^k}.
\label{eq:gamma-k}
\end{equation}

{\bf\punct Averaging operators.}
An element of $M$ having different eigenvalues can be reduced to one of subalgebras $A_k$ by 
a conjugation by some $h\in H$. This element  $a\in A_k$ is defined up to 
the action of $W_k$.

Define a function $\Delta_k$ on $A_k$ by the  Vandermonde expression
$$
\V_k(a^{(k)})=\V(e^{i\phi_1},\dots, e^{i\phi_{n-2k}}, e^{z_k}, e^{-\ov z_k},\dots, e^{z_1}, e^{-\ov z_1}).
$$
Notice that $\V_k(a^{(k)})$ are pure real or pure imaginary depending on $n$, $k$; also  the sign
of $\V_k(a^{(k)})$
depends on the choose of order of eigenvalues.

On the other hand, denote by $B_k\subset H$ the subgroup fixing all elements of 
$A_k$. It consists of block diagonal real matrices with $n-2k$ blocks of size $1\times 1$
and $k$ blocks of size $2\times 2$ having the form 
$\begin{pmatrix}
  c&d\\-d&c
 \end{pmatrix}$, $c^2+d^2\ne 0$.
 Consider a map
 $A_k\times (H/B_k)\to M$ given by 
 \begin{equation}
 (a^{(k)},y)\mapsto ya^{(k)}y^{-1}.
 \label{eq:a-y}
 \end{equation}
 An element $y\in H$ is a matrix  determined up to an equivalence 
 $y\sim yb$, where $b\in B_k$. Therefore the map is well-defined. 
 Equip $H/B_k$ with an invariant measure $d_ky$, which will be normalized several lines below.
 For any function $f\in C^\infty_c(G)$ we define a collection of functions $ I_kf$ on $A_k$
 (integrals of $f$ over orbits of $H$):
 $$
 I_kf(a^{(k)})=\int_{H/B_k} f\bigl(y a^{(k)}y^{-1}\bigr)\,d_ky
 $$
 (functions $ I_kf$ are well-defined on the set $\Delta_k(a^{(k)}\ne 0$).
 Then under a certain normalization of measures $dx$ and $d_k y$
 we have the identity
\begin{equation}
\int_{G/H} f(x)\, dx=
\sum_{k=0}^{[n/2]} \gamma_k \int_{A_k} I_k f(a^{(k)})\, |\V(a^{(k)})|^2\,d a^{(k)}.
\label{eq:weyl-inegration}
\end{equation}
This is a version of the Weyl integration formula. At each point $a^{(k)}\in A_k$ the volume form
on $G/H$
is a product of forms
$$
|\V(a^{(k)})|^2\,da^{(k)} \wedge d_k y
$$
and this explains a normalization of $d_k y$. The factors $\gamma_k=1/{\# W_k}$
arise because the map (\ref{eq:a-y}) covers each point of its image $\# W_k$ times.
See \cite{San}, Sect.4.

\sm

{\bf \punct Normalized averaging operators.}
For each $k$ we define the expression $\epsilon_k(a^{(k)})$ on $A_k$ by
$$
\epsilon_k(a^{(k)}):=\prod_{1\le p<q\le {n-2k}} \sgn\sin\Bigl(\frac{\phi_p-\phi_q}2\Bigr).
$$

For each $k$ we define an {\it averaging operator}
from $C_c^\infty(G)$ to the space of functions on $A_k$ by
\begin{equation}
\Xi_k f(a^{(k)}):= \epsilon_k(a^{(k)})\bigl(\det a^{(k)}\bigr)^{-(n-1)/2} \ov{\V_k(a^{(k)})}
 \int_{H/B_k} f(y a^{(k)}y^{-1})\,d_ky.
\label{eq:average-operator}
\end{equation}


The functions $\Xi_k f$ are invariant with respect to the action of the subgroups 
$S_k$, $S_{n-2k}\subset W_k$ and  change signs under the reflections $R_m$, see
(\ref{eq:Rm}).
The  hypersurfaces $t_m=0$ divide $A_k$ into a union of 
$2^k$
'{\it octants}'
\begin{equation}
\Bigl\{
t_1\lessgtr 0,\quad\dots\quad t_k\lessgtr 0.
\label{eq:octant}
\end{equation}
 Generally, $\Xi_kf(a^{(k)})$
is discontinuous on these hypersurfaces. A function  $\Xi_k f(a^{(k)})$ admits a 
smooth continuation to the closure of each 'octant' and the operator $\Xi_k$
is a bounded operator in a natural sense (\cite{War}, Corollary 8.5.1.2).

An intersection $A_k\cap A_{k+1}$ of Cartan subspaces  is a hypersurface 
in each subspace and $\Xi_k f$, $\Xi_{k+1} f$ satisfy some boundary conditions on this hypersurface
(see Subsect. \ref{ss:glueing}).

\sm

Thus we get an operator
$\Xi$, which sends a function $f\in C_c^\infty(G/H)$
to a collection of functions
$$
\Xi:f\mapsto
(\Xi_0 f, \dots, \Xi_{[n/2]} f).
$$
An $H$-invariant distribution on $G/H$ determines a linear functional
on the image of $\Xi$.

\sm

{\bf \punct Differential Vandermonde.}
Consider the following $n$-tuple of first order differential operators
on $A_k$:
\begin{multline}
 \frac\partial{i \partial \phi_1},\dots,  \frac\partial{i \partial \phi_{n-2k}}, \,
 \frac 12 \Bigl(\frac \partial{\partial t_k}+ \frac\partial{i \partial \theta_k} \Bigr),\, 
 \frac 12 \Bigl(-\frac \partial{\partial t_k}+ \frac\partial{i \partial \theta_k} \Bigr),
 \dots, \\
  \frac 12 \Bigl(\frac \partial{\partial t_2}+ \frac\partial{i \partial \theta_2} \Bigr), \,
 \frac 12 \Bigl(-\frac \partial{\partial t_1}+ \frac\partial{i \partial \theta_1} \Bigr).
 \label{eq:X-n}
\end{multline}
Denote them by $X_1$, \dots, $X_n$.
We define the differential operator
$$
\V(\partial)=\V_k(\partial):=\V(X_1,\dots, X_n).
$$
Any function 
$F=\V_k(\partial)\,\Xi_k f(a^{(k)})$
on $A_k$ satisfies the following properties (see \cite{San}, (8.2)):

\sm

A$^\circ$. $F$ is skew-symmetric with respect
to the subgroup $S_{n-2k}\subset W_k$ and invariant with respect to 
of $S_k\ltimes \Z_2^{k}$.

\sm

B$^\circ$. 
$F$ admits a  continuous extension to  the whole $A_k$;

\sm

C$^\circ$. $F$ is smooth on the closure of each 'octant' (\ref{eq:octant}).

\sm

There is a constant $\gamma_*$ depending on a normalization of the measure on $G/H$ such that
\begin{equation}
\V_{[n/2]} \Xi_{[n/2]} f(0) =\gamma_* f(e)
\label{eq:fund}
\end{equation}
for any $f\in C_c^\infty(G/H)$.


\sm

{\bf\punct Distributions $\Lambda_p$.}
Recall that a function $\ctg \phi/2$ determines a distribution on
the circle $\R/2\pi \Z$ by the formula
$$
f\mapsto \mathrm{p.v.}\int_{-\pi}^\pi \ctg\frac \phi 2\, f(\phi)\, d\phi.
$$

Consider an even-dimensional  torus $\T^{2m}$ with standard coordinates
$e^{i\phi_1}$, \dots, $e^{i\phi_{2m}}$.
Define the distribution $\Lambda_{2m}$ by
\begin{equation}
\Lambda_{2m}(\phi):= \frac{(-i)^m }{(2\pi)^{m} 2^m m! } 
\sum_{\sigma\in S_{2m}}(-1)^\sigma \prod_{j=1}^{m}
\Bigl(\ctg\frac {\phi_{\sigma(2j-1)}}2 \cdot 
\delta\bigl(\phi_{\sigma(2j-1)}+\phi_{\sigma(2j)}\bigr)\Bigr).
\label{eq:Lambda-even}
\end{equation}
For an odd-dimensional torus $\T^{2m+1}$ we define the distribution $\Lambda_{2m+1}$ by
\begin{multline}
\Lambda_{2m+1}(\phi):=\frac{(-i)^m}{(2\pi)^{m} 2^m m!} 
\sum_{\sigma\in S_{2m+1}}(-1)^\sigma \times \\\times
\delta\bigl(\phi_{\sigma(2m+1)}\bigr)\cdot 
\prod_{j=1}^{m}
\Bigl(\ctg\frac {\phi_{\sigma(2j-1)}}2 \cdot 
\delta\bigl(\phi_{\sigma(2j-1)}+\phi_{\sigma(2j)}\bigr)\Bigr).
\label{eq:Lambda-odd}
\end{multline}

{\bf\punct  Formulas for projectors.\label{ss:th}}
The  purpose of the present paper is the following formula.

\begin{theorem}
	\label{th}
Invariant distributions $\Theta_r:C_c^\infty(G/H)\to \C$  determining the projectors $\Pi_0$, \dots, $\Pi_{[n/2]}$
in {\rm (\ref{eq:LLL})--(\ref{eq:PPP})} are given by the formula
\begin{multline}
\gamma_*\la f,\Theta_r\ra_G=\frac{(-1)^{n(n-1)/2}  (n-2r)!}{([n/2]-r)!}\sum_{k=0}^r
(-1)^k 4^{r-k}
\times\\ \times
\gamma_k
\Bigl\la \V_k(\partial) \Xi_k f(\phi,\theta,t),\,\, \Lambda_{n-2k}(\phi)\cdot
 \prod_{j=1}^k \delta(t_j)\delta(\theta_j)
\Bigr\ra_{A_k}.
\end{multline}
\end{theorem}

In particular, the projector corresponding to the discrete series is determined by the distribution
$$
\gamma_*\Theta_0(f)=\frac {(-1)^{n(n-1)/2}(n-2r)!}{[n/2]!}\bigl\la \V_0 \Xi_0 f, \Lambda_n \bigr\ra_{A_0}.
$$

{\bf\punct Remarks on the general problem for pseudo-Riemannian symmetric spaces.}
The problem of separation of spectra is reduced to integration of spherical 
distributions as functions of parameters with respect to the Plancherel measure.

The spaces $\GL(n,\C)/\GL(n,\R)$ are representatives of  spaces
$G(\C)/G(\R)$, where $G(\C)$ is a complex semisimple (reductive) Lie group and $G(\R)$ 
is its real form. The Plancherel formula for such spaces was obtained by
Sh.~Sano, N.~Bopp, and P.~Harink \cite{San2},  \cite{BP}, \cite{Har}. 
The spaces $G(\C)/G(\R)$ are
close relatives and spherical distributions admit elementary expressions
similar to (\ref{eq:kappa-1})-(\ref{eq:kappa-2}). 

In a recent preprint \cite{Ner-U} the similar problem was solved for 
$L^2$ on the pseudo-unitary group $\U(p,q)$ (formulas for projectors have another form but calculations are 
similar).
In the case of real semisimple groups
spherical distributions 
are characters, according Harish-Chandra characters are locally integrable functions
admitting elementary expressions.

For general pseudo-Riemannian semisimple symmetric spaces the way used here and in  \cite{Ner-U}
is impossible.
On the other hand, calculations of V.~F.~Molchanov in rank one case
\cite{Mol2}
 indicate that a general formulation of the problem requires an improvement.

\sm

{\bf \punct Structure of the paper.}
The proof of Theorem \ref{th} is based on the Plancherel formula obtained 
by Sano \cite{San} and also on his proof. For this reason, we must expose
numerous elements of the paper \cite{San}. 
In the Section 2 we establish a skew-symmetric counterpart
of the formula
$$
\sum_{n=-\infty}^\infty e^{in\phi} =2\pi \delta(\phi).
$$
Section 3 contains preliminaries
and  Section 4 evaluation of the projectors.

\sm

{\bf Acknowledgements.} I am grateful to H.Upmeier for discussions of this topic.


\section{A skew-symmetric analog of the delta-function}

\COUNTERS

{\bf \punct  Distributions $\Lambda$.%
\label{th:Lambda}} 
For integers $a_1>\dots >a_p$ denote 
\begin{equation}
\e^a_p=
\e^{a_1,\dots,a_p}_p(e^{i\phi_1}, \dots,e^{i\phi_p}):=
\sum_{\sigma\in S_p} (-1)^\sigma e^{i\sum_m a_{\sigma(m)} \phi_m }
.
\label{eq:e-a}
\end{equation}
Denote by $\Lambda_p$ the following distribution on a torus $\T^p$:
\begin{equation}
\cL_p(e^{i\phi_1}, \dots,e^{i\phi_p}):=\frac 1{(2\pi)^p}\sum_{a_1>\dots> a_p} \e^a_p.
\label{eq:Lambda-again}
\end{equation}

\begin{theorem}
\label{th:2}
We have
 $$\cL_p=
 \Lambda_p,$$
 where $\Lambda_p$ is given by {\rm (\ref{eq:Lambda-even}), (\ref{eq:Lambda-odd})}.
\end{theorem}

 For a proof of this statement 
we need some combinatorial lemmas.

\sm

\begin{figure}
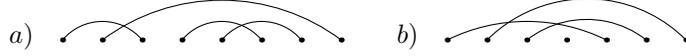

$$a)\quad\epsfbox{diagram.2}\qquad b)\quad \epsfbox{diagram.3}$$
\caption{A matching. Even case and odd case.%
	\label{f:matching}}
\end{figure}

{\bf\punct Matchings.} Consider a $p$-element ordered set $C$, it is convenient
to assume $C=\{p,p-1,\dots,1\}$. If $p$ is even, we say that
a {\it matching} of $C$ is a partition
\begin{equation}
\zeta:\, C=\{c_1,c_2\}\sqcup \{c_3,c_4\}\sqcup \dots
\label{eq:matching-1}
\end{equation}
  of $C$ into two-element subsets. If $p$ is odd, 
then a {\it matching} is a partition of $C$ into $(p-1)/2$ two-element subsets and
one single point subset,
\begin{equation}
\zeta:\, C=\{c_1,c_2\}\sqcup \dots \sqcup \{c_{p-2},c_{p-1}\}\sqcup \{c_p\}
\label{eq:matching-2}
\end{equation}
 We draw matchings as diagrams with arcs, see Figure \ref{f:matching}.
Denote by 
$$\Match(C)=\Match_p$$
the set of all matchings.
Define the standard matching $\zeta_0$ by
$$
\zeta_0=
\begin{cases}
\{p,p-1\}\sqcup\{p-2,p-3\}\sqcup\dots \sqcup \{2,1\},\qquad \text{if $p$ is even}; 
\\
\{p,p-1\}\sqcup\{p-2,p-3\}\sqcup\dots \sqcup \{3,2\}\sqcup \{1\},\qquad \text{if $p$ is odd}.
\end{cases}
$$

\begin{figure}
$$
\epsfbox{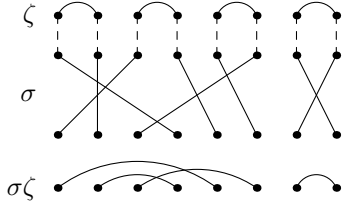}
$$
\caption{Action of symmetric group on $\Match_p$.\label{f:action}}
\end{figure}

 The symmetric group $S_p$ acts on $\Match_p$ 
in a natural way. It is convenient to imagine 
this action as on Figure \ref{f:action} (we glue a diagram of matching
with a diagram of a permutation).

Consider a matching (\ref{eq:matching-1}) or (\ref{eq:matching-2}). We say that pairs 
$\{c_\alpha, c_{\alpha+1}\}$ and $\{c_\beta, c_{\beta+1}\}$ are {\it interlacing}
if precisely one point $c_\beta$ or $c_{\beta+1}$ lies between
$c_\alpha$ and $c_{\alpha+1}$ (in the sense of the ordering of $C$). On Figure \ref{f:matching}
an interlacing of pairs corresponds to an intersection of arcs.
If $\#C=p$ is odd, we say that a distinguished element $\{c_{p}\}$
of a matching {\it interlaces} a pair $\{c_\alpha, c_{\alpha+1}\}$
if $c_{p}$ lies between $c_\alpha$ and  $c_{\alpha+1}$.

We say that a matching $\zeta$ is {\it even} (respectively, {\it odd})
if the number of pairs of interlacing elements of $\zeta$ is even (respectively, odd).
Denote the parity by $(-1)^\zeta$.

\begin{lemma}
	For any $p$,
	\begin{equation}
\Sigma(p):=	\sum_{\zeta\in\Match_p} (-1)^\zeta=1.
\label{eq:1111}
	\end{equation}
\end{lemma}

{\sc Proof.} First, let $p$ be even.
 Define an involution $J$ of $\Match_p$ in the following way.
 
 \sm
 
 --- If $\{p,p-1\}$ is not an element of a matching $\zeta$, we apply to $\zeta$ the transposition
 $(p,p-1)\in S_p$.
 
 \sm
 
 --- If $\{p,p-1\}$ is an element of $\zeta$ and $\{p-2,p-3\}$ not, we apply to $\zeta$ the transposition
 $(p-2,p-3)\in S_p$. Etc.
 
 \sm
 
 --- $J$ fixes the matching $\zeta_0$.
 
 \sm
 
 The involution $J$ changes parity of all matchings $\zeta\ne \zeta_0$. This implies our statement.
 
 Next, let $p$ be odd. Forgetting a singleton $\{c_p\}$ in $\zeta\in\Match_p$, we get a map 
 $h:\Match_p\to\Match_{p-1}$.
  Moreover,
 $$
 (-1)^{\zeta}=(-1)^{c_p+1} (-1)^{h(\zeta)}
 .
 $$
For $c_p=1,3,5,\dots, p$, we have $ (-1)^{\zeta}=(-1)^{h(\zeta)}$,
for $c_p=2,4,\dots, 2p-1$, the signs are different.
Therefore $\Sigma(p)=\Sigma(p-1)$.
\hfill $\square$

\sm

\begin{lemma}
	\label{l:sigma-pi}
For any permutation $\sigma\in S_p$
$$
(-1)^{\sigma \zeta_0}=(-1)^\sigma\cdot \sgn \bigl(\sigma(p)-\sigma(p-1)\bigr)
\cdot
\sgn \bigl(\sigma(p-2)-\sigma(p-3)\bigr)\dots
$$
\end{lemma}

{\sc Proof.} We imagine  $\sigma\in S_p$ as a bipartite graph as on Figure \ref{f:action}.
Inversions in permutations correspond to intersections of arcs. Parity of a permutation
is parity of number of intersections. A matching $\sigma \zeta_0$ is obtained by gluing 
the matching $\zeta_0$ and the permutation $\sigma$. Multiplication of $(-1)^\sigma$
by $\prod_{j=0}^{p-1} \sgn \bigl(\sigma(p-2j)-\sigma(p-2j-1)\bigr)$ means that we do not 
take to account possible inversions in pairs $(p-2j, p-2j-1)$. Now the statement
must be clear from Figure \ref{f:action}.
\hfill $\square$

\sm

{\bf \punct A transformation of the expression for $\Lambda_p$.}
Now we can write formulas (\ref{eq:Lambda-even}), (\ref{eq:Lambda-odd})
in the form
\begin{multline}
\Lambda_p(e^{i\phi_1},\dots, e^{i\phi_p})= (-i)^{[p/2]}
\sum_{\zeta\in \Match_p}
(-1)^\zeta \prod_{\{\alpha,\beta\}\in \zeta} \cot (\phi_\alpha/2)
\bigl(\delta(\phi_\alpha)+\delta(\phi_\beta)\bigr) 
\times\\
\times
\prod_{\{\gamma\}\in \zeta} \delta(\phi_\gamma)
.
\label{eq:Lambda-new}
\end{multline}
The  product $\prod_{\{\gamma\}\in \zeta}$ is taken over the set
 consisting of 0 or 1 elements, i.e., for
even $p$ the product equals 1 and for odd $p$ it consists of 1 factor.

\sm

{\bf \punct Proof of Theorem \ref{th:2}.}

\begin{lemma}
	Let $a_1>\dots> a_{p}$, and $\e^a_p$ be given by {\rm(\ref{eq:Phi})}.
	
{\rm a)} $p=2m$ be even.	 Then
	\begin{equation}
i^m\Bigl\la \e^a_p,\,\, \prod_{l=1}^m
 \Bigl(\ctg\frac {\phi_{2l-1}}2\cdot \delta(\phi_{2l-1}+\phi_{2l})\Bigr)
	\Bigr\ra_{\T^{2m}}=(2\pi)^{2m} 2^m m!	.
	\label{eq:1}
	\end{equation}

{\rm b)}
$p=2m+1$ be odd.	 Then
\begin{equation}
i^m \Bigl\la \e^a_p,\,\, \delta(\phi_{2m+1}) \prod_{l=1}^m
\Bigl(\ctg\frac {\phi_{2l-1}}2\cdot \delta(\phi_{2l-1}+\phi_{2l})\Bigr)
\Bigr\ra_{\T^{2m}}=(2\pi)^{2m+1} 2^m m!.
\label{eq:2}
\end{equation}	
	
\end{lemma}

{\sc Proof.}
We prove a), a proof of b) is  same.
Recall that
$$
i\ctg \phi/2=\sum_{n>0}(e^{-in\phi} -e^{in\phi} ).
$$
Therefore
\begin{equation}
i^m \Bigl\la (-1)^\sigma e^{i a_{\sigma(1)}\phi_1+i a_{\sigma(2)}\phi_2+\dots  },\,
\prod_{l=1}^m
\Bigl(\ctg\frac {\phi_{2l-1}}2\cdot \delta(\phi_{2l-1}+\phi_{2l})\Bigr)
\Bigr\ra =
\label{eq:predlast}
\end{equation}
\begin{multline}
\qquad
=
i^m (2\pi)^m
(-1)^\sigma 
\Bigl\la
e^{i(a_1-a_2)\phi_1+i(a_3-a_4)\phi_3+\dots}, \prod_{l=1}^m \ctg\frac {\phi_{2l-1}}2
\Bigr\ra
=\\
= (2\pi)^{2m} (-1)^\sigma  \prod_{l=1}^m \sgn (a_{\sigma(2l-1)}-a_{\sigma(2l)})
= (-1)^{\sigma\zeta_0} (2\pi)^{2m}
.
\label{eq:last1}
\end{multline}
Thus (\ref{eq:last1}) is
$$
(2\pi)^{2m}
\sum_{\sigma\in S _p} (-1)^{\sigma \zeta_0}=2^m m! (2\pi)^{2m} \sum_{\zeta\in \Match_p} (-1)^{\zeta}
.
$$
The latter sum was evaluated in Lemma \ref{l:sigma-pi}.
\hfill $\square$

\sm

To be definite, assume that $p$ is even.
Since $\e^a_p$ is skew-symmetric, replacing  
$$I(\phi):=\prod
 \ctg\frac {\phi_{2l-1}}2\cdot \delta(\phi_{2l-1}+\phi_{2l})$$
 in (\ref{eq:1})
 by its skew-symmetrization 
 $$
 J:=\sum_{\sigma\in S_p} (-1)^\sigma I(\phi_{\sigma(1)},\dots, \phi_{\sigma(p)})
 $$
 we get the same right-hand side
 multiplied by $p!$. Next, we replace $\e^a_p$ by $\exp(\sum i a_k\phi_k)$,
 after this the right hand side is divided by $p!$. Thus we get Fourier coefficients of
 the complex conjugate distribution $\ov J$.
This is
Theorem \ref{th:2}.

\section{The Plancherel formula}

\COUNTERS

This section contains preliminaries from  \cite{San}.

\sm

{\bf \punct Choose of a sign in the average operator.%
	\label{ss:sign}}
For even $n$
the right hand side of formula (\ref{eq:average-operator}) is defined up to a sign.
Let us fix it.
  Denote 
$$
\nu_k(a^{(k)})=
\prod_{1\le p<q\le {n-2k}} \sin\Bigl(\frac{\phi_p-\phi_q}2\Bigr).
$$
We have
\begin{multline*}
\nu_k(a^{(k)}) \bigl(\det a^{(k)}\bigr)^{-(n-1)/2}=\prod_{1\le  p<q\le n-2k}
(2i)^{-1}
e^{-i(\phi_p-\phi_q)/2}\bigl(e^{i(\phi_p-\phi_q)}-1 \bigr)\times
\\ \times 
\prod_{p=1}^{n-2k} e^{i\phi_p(n-1)/2} \cdot 
\prod_{r=1}^k e^{i\theta_r(n-1)}
.
\end{multline*}
Collecting factors $e^{i\phi_m}$ we get a  continuous expression
$$
e^{-i(\phi_2+2\phi_3+\dots+(n-2k-1)\phi_{n-2k})} \prod e^{i\theta_r(n-1)} \prod 
\bigl(e^{i(\phi_p-\phi_q)}-1 \bigr)
,
$$
which is well-defined up to global change of a sign.
We fix its sign assuming that $\nu_k(a^{(k)})$ is positive
if the order $e^{i\phi_1}$, \dots, $e^{i\phi_{n-2k}}$
are located clock-wise.
It remains to set
 $$
 \epsilon_k(a^{(k)})\bigl(\det a^{(k)}\bigr)^{-(n-1)/2}:= \frac{\nu_k(a^{(k)})
 	 \bigl(\det a^{(k)}\bigr)^{-(n-1)/2}}{|\nu_k(a^{(k)})|}
 .
 $$

{\bf\punct Conditions of glueing.%
	\label{ss:glueing}}
Consider a Cartan subspace $A_k$ (see (\ref{eq:A-k}) with the standard
coordinates 
\begin{equation}
\text{$\phi_1$, \dots, $\phi_{n-2k-2}$, $\phi_{n-2k-1}$, $\phi_{n-2k}$,
	 $t_k$, $\theta_k$, \dots, $t_1$, $\theta_1$}
\label{eq:coordinates-1}
\end{equation}
and the Cartan subspace $A_{k+1}$
with coordinates
\begin{equation}
\text{$\phi_1$, \dots, $\phi_{n-2k-2}$,  $t_{k+1}$, $\theta_{k+1}$, $t_k$, $\theta_k$,
	 \dots, $t_1$, $\theta_1$.}
\label{eq:coordinates-2}
\end{equation}
The intersection $A_k\cap A_{k+1}$ in  $A_{k+1}$  is given by the equation $t_{k+1}=0$; in
$A_k$ it is defined  by the equation $\phi_{n-2k-1}=\phi_{n-2k}$. Let us change coordinates in $A_k$ and 
in $A_{k+1}$ in the following way. In the both cases we leave coordinates
\begin{equation}
\text{$\phi_1$, \dots, $\phi_{n-2k-3}$,  $t_{k}$, $\theta_{k}$, \dots, $t_1$, $\theta_1$.}
\label{eq:coordinates-3}
\end{equation}
being the same. Also:

\sm

-- in $A_{k+1}$, we rename $t:=t_{k+1}$, $\theta:=\theta_{k+1}$;

\sm

-- in $A_k$ we set $\theta:=(\phi_{n-2k-1}+\phi_{n-2k})/2$, $\tau:=\phi_{n-2k-1}-\phi_{n-2k}$.

\sm

Take a point $b$ of the hypersurface $t=0$ such that other $t_j\ne 0$.
Then for any $N>0$   there are smooth functions $u_j(\cdot)$ depending
on coordinates (\ref{eq:coordinates-3})  and $\theta$ such that in a small
neighborhood of 
$b$ we have asymptotic expansions of the form
\begin{align*}
\Xi_{k+1} f (t,\dots)&=\sum_{j=0}^N u_j t^j+o(t^N), \qquad \text{for $t>0$};
\\
\Xi_{k+1} f (t,\dots)&=\sum_{j=0}^N (-1)^{j+1} u_j t^j+o(t^N), \qquad \text{for $t<0$}; 
\\
\Xi_k f(\tau,\dots)&=i \sum_{m=0}^{[N/2]} (-1)^m u_{2m} \tau^{2m}+o(\tau^N),
\end{align*}
see \cite{San}, Theorem 4.4.

\sm

{\sc Remark.} For $t\ge 0$ denote
$$\Xi^+_{k+1} f (t,\dots)=(\Xi_{k+1} f (t,\dots)+
\Xi_{k+1} f (-t,\dots))/2$$
Then the function
$$
 \wt \Xi f(s):= 
\begin{cases}
\Xi_{k+1}^+ f (\sqrt s,\dots), \qquad \text{for $s\ge 0$};
\\
i\Xi_k f(\sqrt{-s},\dots),\qquad\text {for $s\le 0$}
\end{cases}
$$
is a smooth function near $s=0$. 
Appearance of the factor $i$ is artificial, it is related to the
normalization of the factor $\V(a^{(k)})$ in (\ref{eq:average-operator}).
See also \cite{Ten} for elementary
explanations, our asymptotics can be reduced to the case $p=1$, $q=2$ of that paper.
\hfill $\boxtimes$

\sm

{\bf\punct Spherical distributions}, see \cite{San}, Sect. 7.
Points of the spectrum of $L^2(G/H)$ are enumerated by an integer $r=0$, \dots, $[n/2]$
and signatures of type $r$:
\begin{equation}
(c,l)=
(c_1,\dots,c_{n-2r}, l_1,\ov l_1\dots, l_r,\ov l_r).
\label{eq:signature}
\end{equation}
Here $c_1>c_2>\dots>c_{n-2r}$ are integers, and
$$
l_p=(m_p-i\lambda_p)/2,
$$
where $m_p\in \Z$, $\lambda_p\in \R$, and $\lambda_1>\lambda_2>\dots >\lambda_r>0$.
To write a distribution corresponding to a given signature $(c,l)$, we need some notation.

Let
$\phi_1$, $\phi_2$ be defined modulo $2\pi$.
Then $(\phi_1-\phi_2)/2$ is defined modulo $\pi$.
We will use two ways to define $(\phi_1-\phi_2)/2$ modulo $2\pi$  setting
\begin{equation}
\Bigl\lfloor \frac {\phi_1-\phi_2}2 \Bigr\rfloor\in (-\pi/2,\pi/2),
\qquad \qquad
\Bigl\lceil \frac {\phi_1-\phi_2}2 \Bigr\rceil\in (0,\pi).
\label{eq:floor}
\end{equation}
Next, for $l=(m-i\lambda)/2$ we define a function
\begin{equation}
 \xi(l;e^z):=e^{im\theta}(e^{i\lambda t}- e^{-i\lambda t}), \qquad \text{where $e^z=e^{t+i\theta}$},
 \label{eq:xi}
\end{equation}
and a function $D(l;e^{i\phi_1},e^{i\phi_2})$. If $m\in 2\Z$,
we set
\begin{equation}
 D(l;e^{i\phi_1},e^{i\phi_2}):=2 e^{im (\phi_1+\phi_2)/2} \,
 \frac{\ch \lambda \bigl(\bigl|\lfloor (\phi_1-\phi_2)/2 \rfloor \bigr|-\pi/2) }
 {\sh \pi\lambda/2}.
\end{equation}
If $m\in 2\Z+1$, then (formula (8.6) in \cite{San} contains a typos).
\begin{multline}
 D(l;e^{i\phi_1},e^{i\phi_2}):=-2 e^{i(m-1) (\phi_1+\phi_2)/2} 
  \times\\\times
 e^{i\phi_2+\lceil (\phi_1-\phi_2)/2)\rceil}\,
 \frac{\sh \lambda(\lceil (\phi_1-\phi_2)/2)\rceil-\pi/2\bigr)}
 {\ch \lambda\pi/2}
 .
\end{multline}
Functions $D_l$ are continuous on the torus $\T^2$, smooth outside the diagonal 
$\phi_1=\phi_2$,
and have a kink on the diagonal.

We write a spherical distribution $\Phi^r_{c,l}$ corresponding to a signature
(\ref{eq:signature}) in the form
\begin{equation}
\la f, \Phi^r_{c,l}\ra=  \sum 
\gamma_k \int_{A_k} \Xi_k f(a^{(k)})\, \kappa_k(c,l;a^{(k)})\, da^{(k)}
\label{eq:Phi}
,\end{equation}
where
$\gamma_k$ are given by (\ref{eq:gamma-k})  and 
$\kappa_k(c,l;a^{(k)})$ will be defined now.

\begin{figure}
$$
\epsfbox{diagram2.1}
$$
\caption{A diagram $\frS$.\label{f:diagram}}
\end{figure}

Let $k\le r$. We consider diagrams
$\frS$ of
the form shown on Fig.\ref{f:diagram}.
The upper row consists of $n-2r$ white circles corresponding to the parameters 
$c_1$, \dots, $c_{n-2r}$ and $2r$ black circles corresponding to the parameters
$l_1$, $\ov l_1$, $l_2$, $\ov l_2$, \dots, $l_r$, $\ov l_r$. 
The lower row consists of $n-2k$ white boxes corresponding to the coordinates
$\phi_1$, \dots $\phi_{n-2k}$ and 
$2k$ black boxes corresponding to the coordinates $z_1$, $-\ov z_1$,  $z_2$, $-\ov z_2$,
\dots, $z_k$, $-\ov z_k$.
In the upper row,
horizontal dotted arcs link elements of pairs ($l_q,\ov l_q$);
in the lower row, elements of 
pairs ($z_\gamma, -\ov z_\gamma$).
Vertical arcs establish a bijection between  
circles  and a boxes. This bijections are not arbitrary,
vertical arcs must belong to one of the following 3 types:

\sm

Type 1. An arc connecting  a white circle 
and  a white box. 
We denote such arcs by $[c_p;\phi_\alpha]$ according the variables on its ends.
The number of such arcs is $n-2r$.

\sm

Type 2. Pairs of arcs from  linked black circles  to white boxes,
$[l_p,\phi_\alpha]$, $[\ov l_p, \phi_\beta]$, where $\alpha<\beta$.
The number of such pairs is $r-k$.

\sm

Type 3. Pairs of non-intersecting arcs from linked black circles to linked black boxes
$[l_p,z_\gamma]$, $[\ov l_p,-z_\gamma]$. The number of such pairs is
$k$.

\sm

Denote the set of such diagrams by $\Omega(r,k)$.

Since $\frS$ determines a permutation, it has a well-defined sign $(-1)^\frS$.
This number is also a parity of the number of intersections of 
the vertical arcs, notice that we can do not take to an account arcs of Type 3.

\sm


\sm

Functions $\kappa_k$ determining spherical distributions (see\ref{eq:Phi}) are given by formula: 
\begin{equation}
\kappa_k(c,l;a^{(k)}):=0 \qquad \text{for $k>r$};
\label{eq:kappa-1}
\end{equation}
\begin{multline}
\kappa_k(c,l;a^{(k)})=\sum_{\frS\in\Omega(r,k)}
\prod_{[c_p,\phi_\alpha]\in \frS} e^{i c_p\phi_\alpha}
\prod_{[l_q,z_\gamma]\in \frS}
\xi(l_q; e^{z_\gamma})
\times\\\times
\prod_{[l_q;\phi_\alpha],[\ov l_q, \phi_\beta]\in\frS}
D(l_q;e^{i\phi_\beta}, e^{i\phi_\gamma})
\qquad\text{for $k\le r$},
\label{eq:kappa-2}
\end{multline}
the summation is taken over all diagrams $\frS$, the product is taken over
all pieces of a given diagram $\frS$
(see \cite{San}, equation (7.6), and the expression $\kappa$ at the last row of
\S7).

\sm
 
 {\sc Remark.}
The expressions $\kappa_k(c,l;a^{(k)})$ are invariant
with respect to the subgroups $S_{n-2k}$, $S_k\subset W_k$ acting on $A_k$ and
and changes a sign under reflections (\ref{eq:Rm}).
 Also, they are invariant with respect to permutations of $c_1$, \dots, $c_{n-2r}$
and with respect to permutations of $l_1$, \dots $l_r$. They change a sign
under the reflections
\begin{equation}
(\dots,l_{j-1},l_j,l_{j+1},\dots)\mapsto  (\dots,l_{j-1},\ov l_j,l_{j+1},\dots).
\label{eq:reflection-l}
\end{equation}

\sm

{\bf \punct Functions $\kappa(c,l;a^{(k)})$ as eigenfunctions of symmetric differential operators.}
Notice, that the functions $\xi$, $D$ can be written as
\begin{align}
&\xi(l;z)=e^{zl+\ov z\ov l}- e^{\ov z l+ z\ov l};
\label{eq:xi-alt}
\\
&D(l;e^{i\phi_1}, e^{i\phi_2})
=\frac{2} {e^{2\pi i l}-1}e^{i\phi_1 l+i\phi_2 \ov l}
-\frac{2} {e^{2\pi i \ov l}-1} 
e^{i\phi_1 \ov l+i\phi_2 l},
\label{eq:D-alt}
\end{align}
in the second expression we choose $\phi_1/2>\phi_2/2$.
 

Therefore locally the functions
$\kappa_k(\phi,z)$ are linear combinations of exponents
of linear functions.

Consider a  symmetric polynomial
$S(x_1,\dots,x_n)$. Substituting the first order differential operators 
(\ref{eq:X-n}) to $S(\cdot)$ we come to the
identity
$$
S(X_1,\dots,X_n)\kappa_k(c,l;a^{(k)})=S(\lambda_1,\dots,\lambda_n)\,\kappa_k(c,l;a^{(k)}),
$$
which  holds
outside hypersurfaces $\phi_p=\phi_q$. Moreover, this identity is valid 
in a distributional sense. Precisely, for any $f\in C^\infty(G)$
\begin{multline*}
\sum_{k=0}^{[n/2]}\gamma_k\int_{A_k} S(-X) \Xi_k f(a^{(k)})\cdot \kappa_k(a^{(k)})\,da^{(k)}
=\\=\sum_{k=0}^{[n/2]}\gamma_k\int_{A_k}  \Xi_k f(a^{(k)})\cdot S(X) \kappa_k(a^{(k)})\,da^{(k)}.
\end{multline*}
Boundary terms, which appear after integration by parts,  cancel due the gluing
conditions, see \cite{San}, Lemma 6.2. Sano also establishes a more general integration by parts identity
(it will be used below).
Let $S(x_1,\dots,x_n)$,  $T(x_1,\dots,x_n)$ be homogeneous polynomials,
which are both symmetric or both skew-symmetric in $x_1$, \dots, $x_n$.
Consider differential operators $S(X_1,\dots,X_n)$, $T(X_1,\dots,X_n)$, where
$X_j$ are given by (\ref{eq:X-n}). Then
\begin{multline}
\sum_{k=0}^{[n/2]}\gamma_k\int_{A_k}  \Xi_k f(a^{(k)})\cdot S(X) T(X) \kappa_k(a^{(k)})\,da^{(k)}
=\\=
\sum_{k=0}^{[n/2]}\gamma_k\int_{A_k} S(-X) \Xi_k f(a^{(k)})\cdot T(X) \kappa_k(a^{(k)})\,da^{(k)}.
\label{eq:by-parts}
\end{multline}

\sm

{\bf \punct The Plancherel Theorem.}
Consider the distributions $\Phi^r_{c,l}$ on $G/H$ given by 
(\ref{eq:Phi}). Denote (see formula(\ref{eq:vandermonde}))
$$
\V(c,l):=\V(c_1,\dots,c_{n-2r},l_1,\ov l_1,\dots, l_r,\ov l_r).
$$
Then the following Plancherel  formula holds.
For any $f\in C_c^\infty(G/H)$,
\begin{multline}
\gamma_* f(e)=
\frac 1{(2\pi)^n}\sum_{r=0}^{[n/2]}
\frac{(n-2r)!\, i^{[n/2]-r} (-1)^r}{([n/2]-r)!}
\times \\ \times
\sum_{c_1>\dots>c_{n-2r}}
\sum_{m_1,\dots,m_r}
\int_{\lambda_1>\dots>\lambda_r>0}
\la f, \Phi^r_{c,l}\ra\, \V(c,l)\,d\lambda_1\dots d\lambda_r.
\label{eq:plancherel}
\end{multline}
where $\Phi^r_{c,l}$ are spherical distributions
given be (\ref{eq:Phi}) and
$\gamma_*$ is the same constant as in (\ref{eq:fund}).
See \cite{San}, Theorem 8.7.
The formula assumes that
summands of the exterior sum $\sum_r$  are absolutely convergent (as  integrals over  measures
on spaces $\Z^{n-2r}\times \Z^r\times \R_+^r$).

\section{Calculations}

\COUNTERS

Thus, we must evaluate summands of the exterior sum in (\ref{eq:plancherel}),
i.e. we must find a distribution given by
\begin{multline}
\gamma_*(2\pi)^n\la f,\Theta_r\ra=
\frac{(n-2r)!\, i^{[n/2]-r} (-1)^r}{([n/2]-r)!}
\sum_{c_1>\dots>c_{n-2r}}\,
\sum_{m_1,\dots,m_r}
\times \\\times
\int_{\lambda_1>\dots>\lambda_r>0}
\Bigl(\sum_{k=0}^{r} \gamma_k
\int_{A_k} \Xi_k f(a^{(k)})\cdot \kappa_k(c,l;a^{(k)}) \,d a^{(k)}\Bigr) \V(c,l)
\,d\lambda_1\dots d\lambda_r.
\end{multline}

{\bf \punct Integration by parts.}
Let $l$ be as above, $l=(m-i\lambda/2)$.
Define functions 
$$
\xi'(l;e^z):=-e^{im\theta}(e^{i\lambda \theta}+e^{-i\lambda \theta}).
$$
If $m\in 2\Z$,
we set
\begin{multline*}
D'(l;e^{i\phi_1},e^{i\phi_2}):=2 e^{im (\phi_1+\phi_2)/2} 
\times \\ \times
\frac{\sh \lambda \bigl(\bigl|\lfloor (\phi_1-\phi_2)/2 \rfloor \bigr|-\pi/2) }
{\sh \pi\lambda/2} \cdot \sgn \bigl(\lfloor (\phi_1-\phi_2)/2 \rfloor\bigr)
.
\end{multline*}
If $m\in 2\Z+1$, then 
\begin{multline}
D'(l;e^{i\phi_1},e^{i\phi_2}):=-2 e^{i(m-1) (\phi_1+\phi_2)/2} 
\times\\\times
e^{i\phi_2+ i\lceil (\phi_1-\phi_2)/2)\rceil}
\frac{\ch \lambda(\lceil (\phi_1-\phi_2)/2)\rceil-\pi/2\bigr)}
{\ch \lambda\pi/2}
.
\end{multline}
A function $D(l;e^{i\phi_1},e^{i\phi_2})$ has a jump
on the diagonal 
$e^{i\phi_1}=e^{i\phi_2}$ and $C^\infty$-smooth outside the diagonal.

Following \cite{San}, consider  functions $\kappa^{\,\prime}(c,l;a^{(k)})$
on the union of $A_k$ given by the formula
\begin{equation}
\kappa_k'(c,l;a^{(k)}):=0 \qquad \text{for $k>r$}.
\label{eq:kappa-prime1}
\end{equation}
\begin{multline}
\kappa_k^{\,\prime}(c,l;a^{(k)})=\\=
\sum_{\frS\in \Omega(r,k)}
(-1)^{\frS}
\prod_{[c_p,\phi_\alpha]\in \frS} e^{i c_p\phi_\alpha}
\prod_{[l_q,z_\delta]\in \frS}
\xi'(l_q; e^{z_\delta})
\times \\ \times
\prod_{[l_q;\phi_\alpha], [\ov l_q,\phi_\beta]\in\frS}
D'(l_q;e^{i\phi_\alpha}, e^{i\phi_\beta})
\qquad \text{for $k\le r$}
\label{eq:kappa-prime2}
.\end{multline}

\sm

{\sc Remark.}
These functions are skew-symmetric with respect to the variables 
$\phi_m$ and invariant with respect to the subgroup 
$S_{k}\ltimes \Z_2^k\subset W_k$. They are skew-symmetric with respect to the
parameters $c_m$ and invariant with respect to the permutations of 
$l_\alpha$ and the reflections (\ref{eq:reflection-l})
\hfill $\boxtimes$

\sm
 
We have
\begin{align}
\V_k(\partial) \,\kappa_k(c,l;a^{(k)})=\V(c,l) \,\kappa_k^{\,\prime}(c,l;a^{(k)});
\\
\V_k(\partial)\, \kappa_k^{\,\prime}(c,l;a^{(k)})=\V(c,l)\, \kappa_k(c,l;a^{(k)}).
\end{align}
These identities hold pointwise outside diagonals 
$e^{i\phi_p}=e^{i\phi_q}$, this easily follows from (\ref{eq:xi-alt}).
Moreover, we have an identity:
\begin{multline*}
\sum_{k\le r} \gamma_k \int_{A_k}\Xi_k f(a^{(k)}) \cdot \V(c,l)  \kappa_k (c,l;a^{(k)})\, da^{(k)}=
\\=
(-1)^{n(n-1)/2}
\sum_{k\le r} \gamma_k  \int_{A_k}\V_k(\partial)\Xi_k f(a^{(k)}) \cdot
\kappa^{\,\prime}_k (c,l;a^{(k)})\, da^{(k)}
\end{multline*}
To obtain this, we apply (\ref{eq:by-parts}) with $S$, $T=\V$.

Thus, 
$$
\gamma_*\cdot(2\pi)^n\la f,\Pi_r\ra=
\frac{(n-2r)!\, i^{[n/2]-r} (-1)^r(-1)^{n(n-1)/2}}{([n/2]-r)!}
\sum_{k}\gamma_k U_{r,k}
,$$
where 
\begin{multline}
U_{r,k}=
\sum_{c_1>\dots>c_{n-2r}}
\sum_{m_1,\dots,m_r}
\\
\int\limits_{\lambda_1>\dots>\lambda_r>0}
\int\limits_{A_k} h_k(a^{(k)}) \cdot \kappa^{\,\prime}_k(c,l;a^{(k)})\,d a^{(k)}
\,d\lambda_1\dots d\lambda_r
\label{eq:U-k}
\end{multline}
where $h_k(a^{(k)})=\V_k(\partial)\Xi_k f(a^{(k)})$.

\sm

{\bf \punct Summation of distributions.}

\begin{lemma}
	Let $h_k(a^{(k)})$ be  skew-symmetric with respect to $S_{n-2k}\subset W_k$,
	symmetric with respect to $S_k\ltimes \Z_2^k$ and smooth in the domain
	$t_1\ge 0$, \dots, $t_k\ge 0$. Then the sum given by
{\rm(\ref{eq:U-k})} equals to
$$
\frac{(2\pi)^n  (n-2r)! (-2)^{r-k}(-i)^{[n/2]-r}}
{(n-2r)!}  \Bigl\la h_k, \Lambda_{n-2k}(\phi) \prod_{m=1}^k \delta(\theta_m)\delta (t_m)\Bigr\ra_{A_k}
.
$$	
\end{lemma}

{\sc Notation for a proof.}
For a diagram $\frS\in\Omega(r,k)$
consider a diagram $\wt\frS$ obtained by removing
vertical 
arcs of Type 1, i.e., arcs  $[c_p,\phi_\alpha]$, see Figure \ref{fig:diagram-2}.a.
Next,  we define the diagram
$ \frS^\circ\in \Omega(r,k)$ obtained from $\wt\frS$ by adding a collection of nonintersectig arcs
from white circles
to free black boxes. Denote by $\Omega^\circ(r,k)$ the set of diagrams of this form,
$\Omega^\circ(r,k)\subset \Omega(r,k)$.

For $\frT\in\Omega^\circ(r,k)$ we denote by
$u_1<\dots<u_{n-2k}$, which are connected with white boxes.

Also for $\frT\in \Omega^\circ(r,k)$ we define two parites, $\epsilon_1(\frT)$, $\epsilon_2(\frT)$;
let
$\epsilon_1(\frT)$ be the parity of number of intersections  between  arcs of Type 1 and arcs of Type 2.
Respectively, 
$\epsilon_2(\frT)$ is the parity of the number of intersections of arcs of Type 2,
$$
(-1)^\frT=(-1)^{\epsilon_1(\frT)} (-1)^{\epsilon_2(\frT)}.
$$

Finally, for $\frT\in \Omega^\circ(r,k)$ we consider 
the diagram $\frT^\square$ obtained by the following operation:
we forget black circles and black boxes, we forget arcs between black circles and black boxes;
we transform any pairs of arcs $[l_q,\phi_\alpha]$, $[\ov l_q,\phi_\beta]$ to an arc
between $\phi_\alpha$, $\phi_\beta$. We denote the set of such diagrams by
$\Omega^\square(r,k)$.

In the same way, for $Q\in \Omega(r,k)$ we define two numbers
$(-1)^{\epsilon_1(\frQ)}$, $(-1)^{\epsilon_2(\frQ)}$ and the collection
$u_1<\dots<u_{n-2k}$.

\sm


\begin{figure}
\qquad\qquad
$
\epsfbox{diagram2.2}
$

a) The diagram $\wt \frS$ for $\frS$ given be Fig.\ref{f:diagram}.

\bigskip

\qquad\qquad
$
\epsfbox{diagram2.3}
$

b) The diagram $\frT=\frS^\circ\in\Omega^\circ(r,k)$

\bigskip

\qquad\qquad
$
\epsfbox{diagram2.4}
$

c) The diagram $\frQ=\frT^\square\in\Omega^\square(r,k)$.
\caption{.\label{fig:diagram-2}}
\end{figure}

{\sc Proof.}
By the symmetry of the expressions $\kappa^{\,\prime}_k$ in the parameters  with respect to
$S_r\ltimes \Z_2^r$  we can 
represent this as
\begin{multline}
U_{r,k}=\\
\frac 1{r!\,2^r}
\sum_{c_1>\dots>c_{n-2r}}
\sum_{m_1,\dots,m_r}\,\,
\int\limits_{\lambda_1,\dots,\lambda_r\in\R}
\int_{A_k} h_k(a^{(k)})
 \kappa^{\,\prime}_k(c,l;a^{(k)})
 \,da_k
\,d\lambda_1\dots d\lambda_r
=\\=
\frac 1{r!\,2^r}
\sum_{\frS\in \Omega(r,k)}
(-1)^{\frS}
\sum_{c_1>\dots>c_{n-2r}}
\sum_{m_1,\dots,m_r}\,\,
\int\limits_{\lambda_1,\dots,\lambda_r\in\R}
\int_{A_k} h_k(a^{(k)})
\times\\\times
\prod_{[c_p,\phi_\alpha]\in \frS} e^{i c_p\phi_\alpha}
\prod_{[l_q,z_\delta]\in \frS}
\xi'(l_q; e^{z_\delta})
\times \\ \times
\prod_{[l_q;\phi_\alpha], [\ov l_q,\phi_\beta]\in\frS}
D'(l_\gamma;e^{i\phi_\alpha}, e^{i\phi_\beta})
\,da_k
\,d\lambda_1\dots d\lambda_r
\label{eq:long}
.\end{multline}



We will use the following  identities for distributions.
Obviously,
\begin{equation}
\sum_m \int_{\lambda>0} \xi'(l;z)\,d\lambda=-(2\pi)^2 \delta(t)\delta(\theta)
.
\label{eq:id1}
\end{equation}
By \cite{San}, Lemma 8.4,
we have
\begin{multline}
\sum_m \int_{\lambda>0}D'(l; e^{i\phi_1}, e^{i\phi_1})  \,d\lambda
=4i\sum_{a_1>a_2} \bigl(e^{ia_1\phi_1+ia_2 \phi_2 }- e^{ia_2\phi_1+ia_1 \phi_2 }\bigr)=
\\=
4i\Lambda_2(e^{i\phi_1},e^{i\phi_2})= (-8\pi) \cot (\phi_1/2) \delta(\phi_1+\phi_2)
.
\label{eq:id2}
\end{multline}

According \cite{San}, we can change order of summation and integration (\ref{eq:long})   in an arbitrary way.
We successively integrate with respect to the following groups
of variables:

\sm

$\bullet$  for each $[l_q,z_\delta]\in \frS$ we integrate
$$
\int d\theta_\delta \int dt_\delta \sum_{m_q} \int d\lambda_q \bigl(\dots\bigr)
$$
applying (\ref{eq:id1});

\sm

$\bullet$  for each pair
$[l_q;\phi_\alpha]$, $[\ov l_q,\phi_\beta]\in\frS$
we integrate
$$
\int d\phi_{\alpha} \int d\phi_{\beta}
 \sum_{m_q} \int d\lambda_q \bigl(\dots\bigr)
$$
applying (\ref{eq:id1}).

\sm

We come to
\begin{multline*}
U_{r,k}=
\frac {(2\pi)^{n-2r} (8\pi^2)^k (-8\pi)^{r-k} } {r!\,2^r}  
\sum_{\frT\in \Omega^\circ(r,k)}
\times\\\times
\biggl\la h_k,
 (-1)^{\epsilon_1(\frT)}
\Lambda_{n-2r}(e^{i\phi_{u_1}},\dots, e^{i\phi_{u_{n-2r}}})
\times \\ \times
(-1)^{\epsilon_2(\frT)}
\prod_{[l_q,z_\gamma]\in \frT}
\delta(t_\gamma) \delta(\theta_\gamma)
\prod_{[l_q;\phi_\alpha], [\ov l_q,\phi_\beta]\in\frT}
\cot(\phi_\alpha/2)\delta(\phi_\alpha+\phi_\beta)
\biggr\ra
.
\end{multline*}
A summand corresponding to $\frT$ actually depends  on $\frT^\square$,
and there are $r!$ different $\frT^\square$ for a given $\frT$.
We also apply the expression (\ref{eq:Lambda-new})
for $\Lambda$ and get
\begin{multline*}
U_{r,k}=
\frac { (2\pi)^{n-2r} (8\pi^2)^k (-8\pi)^{r-k} r!\, (-i)^{[(n-2r)/2]} } {r!\,2^r}  \times
\\ \times
\Biggl\la h_k,\,\,
\prod_{j=1}^k \delta(t_j) \delta(\theta_j) 
\times\\\times
\biggl\{\sum_{\frQ\in \Omega^\square(r,k)}\!\!\!\!\!\!
(-1)^{\epsilon_1(\frQ)} \!\!\!\!\!\!
\sum_{\zeta\in \Match(\{u_1,\dots, u_{n-2r}\}} \!\!\!\!\!\!
(-1)^\zeta \prod_{\{\alpha,\beta\}\in \zeta} \cot (\phi_\alpha/2)
\bigl(\delta(\phi_\alpha+\phi_\beta)\bigr)
\\
\times (-1)^{\epsilon_2(\frQ)}
\prod_{\{\gamma\}\in \zeta} \delta(\phi_\gamma)\cdot
 (-1)^{\epsilon_2(\frQ)}
\prod_{[\beta,\gamma]\in\frQ}
\cot(\phi_\alpha/2)\delta(\phi_\alpha+\phi_\beta)
\biggr\}\Biggr\ra.
\end{multline*}
The sum in the big curly brackets is $\Lambda_{n-2k}\bigl(e^{i\phi_1},\dots, e^{i\phi_{n-2k}}\bigr)$.
This completes the calculation.
\hfill $\square$

\tt

\noindent
Math. Dept., University of Vienna; \\
Institute for Theoretical and Experimental Physics (Moscow) \\
MechMath Dept., Moscow State University\\
Institute for Information Transmission Problems\\
URL: http://mat.univie.ac.at/$\sim$neretin/

\end{document}